\documentclass[12pt]{article}
\usepackage{amsmath,amssymb}
\usepackage{graphicx,psfrag,epsf}
\usepackage{enumerate}
\usepackage{natbib}
\usepackage{url} 

 \newcommand{\Ex}{{\mathbb E}}
 \renewcommand{\Pr}{{\mathbb P}}

\newcommand{\blind}{0}

\begin{document}

\def\spacingset#1{\renewcommand{\baselinestretch}%
{#1}\small\normalsize} \spacingset{1}


\if0\blind
{
  \title{\bf A Prediction Tournament Paradox}
  \author{David J. Aldous\thanks{
    The author gratefully acknowledges support by N.S.F. Grant DMS-1504802}\hspace{.2cm}\\
    Department of Statistics, U.C. Berkeley CA 94720-3860}

  \maketitle
} \fi

\if1\blind
{
  \bigskip
  \bigskip
  \bigskip
  \begin{center}
    {\LARGE\bf A Prediction Tournament Paradox}
\end{center}
  \medskip
} \fi

\bigskip
\begin{abstract}
 In a prediction tournament, contestants ``forecast" by asserting a numerical probability for each of (say) 100  future real-world events.
 The scoring system is designed so that (regardless of the unknown true probabilities) more accurate forecasters  will  likely score better. 
 This is true for one-on-one comparisons between contestants. 
 But consider a realistic-size tournament with many contestants, with a range of accuracies.
 It may seem self-evident that the winner will likely be one of the most accurate forecasters.
 But, in the setting where the range extends to very accurate forecasters,  simulations show this is mathematically false, within a somewhat plausible model.
 Even outside that setting the winner is less likely than intuition suggests to be one of the handful of best forecasters. 
 Though implicit in recent technical papers, 
 this paradox has apparently not been explicitly pointed out before, though is easily explained.  
 It perhaps has implications for the ongoing IARPA-sponsored research programs involving forecasting.
\end{abstract}

\noindent%
\vfill

\newpage
\spacingset{1.00} 
\section{Introduction}
\label{sec:intro}

General non-mathematical background to the topic is best found in the persuasive essay by \cite{essay}.
Study of results of prediction tournaments in recent years has led to the popular book by \cite{fox} and substantial academic literature  -- 
for instance \cite{psych} has 117 Google Scholar citations.
That literature involves serious statistical analysis, but is focussed on the psychology of individual and team-based decision making and the effectiveness of training methods.
This article considers some (quite elementary)  mathematical questions.

In mathematical terms, 
a prediction tournament consists of a collection of $n$ questions of the form ``state a probability for a specified real-world event
happening before a specified date".
In actual tournaments one can update probabilities as time passes, but for simplicity we  consider only a single probability prediction for each  question, and only binary outcomes.
Scoring is by squared error\footnote{In fact tournaments use {\em Brier score}, which is just $2 \times$ the squared error, with modifications for multiple-choice questions.}: if you state probability $q$ then on that question

score = $(1-q)^2$ if event happens;  score = $q^2$ if not. 

\noindent
Your tournament score is the sum of scores on each question.
As in golf one seeks a {\em low} score. 
Also as in golf, in a {\em tournament} all contestants address the same questions; it is not a single-elimination tournament as in tennis.

The use of squared-error scoring is designed so that (under your own belief) the expectation of your score is minimized by stating 
your actual probability belief.
So in the long run  it is best to be ``honest" in that way.

In more detail, with unknown true probabilities $(p_i)$,  if you announce probabilities $(q_i)$
then (see (\ref{Sdef}) below) the true expectation of your score equals
\[ \sum_i p_i(1-p_i) + \sum_i (q_i - p_i)^2 . \]
The first term is the same for all contestants, so if $S$ and $\hat{S}$ are the tournament scores for you and another 
contestant, then
\[
n^{-1/2} (\Ex S - \Ex \hat{S}) = \sigma^2 - \hat{\sigma}^2 \]
where
\[ \sigma : = \sqrt{n^{-1} \sum_i (q_i - p_i)^2 } \]
is your RMS error in predicting probabilities and $\hat{\sigma}$ is the  other contestant's RMS error.
Thus by looking at differences in scores one can, in the long run, estimate relative abilities at prediction,
as measured by RMS error of predicted probabilities.

To re-emphasize, when we talk about prediction ability we mean the ability to
estimate {\em probabilities} accurately; we are not talking about predicting Yes/No outcomes and counting the number of successes, which is an extremely inefficient procedure for comparing prediction ability.

\subsection{The elementary mathematics}
Let us quickly write down the relevant elementary mathematics.
Write $X$ for your score on a question when the true probability is $p$ and you predict $q$:
\[ \Pr(X = (1-q)^2) = p, \quad \Pr(X = q^2) = 1-p . \]
\[ 
\Ex X = p(1-q)^2 + (1-p)q^2 = p(1-p) + (q-p)^2 .
\]
So writing $S$ for your ``tournament score" when the true probabilities of the $n$ events are $(p_i, 1 \le i \le n)$ and you predict $(q_i, 1 \le i \le n)$,
\begin{equation}
\Ex S = \sum_i p_i(1-p_i) +  n \sigma^2 
\label{Sdef}
\end{equation}
where 
\[ \sigma^2 :=  n^{-1} \sum_i (q_i - p_i)^2 \]
is your MSE (mean squared error) in assessing the probabilities. 
Let us spell out some of the implications of this simple formula. 
\begin{itemize}
\item  The first term in (\ref{Sdef}) is the same for all contestants: one could call it the contribution from  ``irreducible randomness". 
\item The formula shows that a convenient way to measure forecast accuracy is via $\sigma$, the RMS (root-mean-square) error of a candidate's forecasts. 
\item The actual score $S$ is random:
\begin{equation}
 S =  \sum_i p_i(1-p_i) +  \sigma^2  + (\mbox{ chance variation })
 \label{3-part}
\end{equation}
where the ``chance variation" has expectation zero. 
Given the scores $S$ and  $\hat{S}$ for you and another contestant, one could attempt a formal test of significance of the hypothesis 
$\sigma < \hat{\sigma}$ that you are a more accurate forecaster.
But making a valid test is quite complicated, because the ``chance variations" are highly correlated. 
\end{itemize}

\subsection{But one-on-one comparisons may be misleading}
In the model and parameters we will describe in section \ref{sec:whowins}, a contestant in a 100-question tournament who is 5\% more accurate than another 
(that is, RMS prediction errors 10\% versus 15\%, or 20\% versus 25\%) will have around a 75\% chance
to score better (and around 90\% chance if 10\% more accurate).  
This is  unremarkable; it is just like the well-known Bradley-Terry style models (see e.g. \cite{hunter}) for sports, where the probability 
A beats B is a specified function of the difference in strengths. 
In the sports setting it seems self-evident that in any reasonable league season or tournament play, the overall winner 
is likely to be one of the strongest teams.  
The purpose of this paper is to observe, in the next section,  that (within a simple model )
this ``self-evident" feature is just plain false for prediction tournaments. 
So in this respect,  prediction tournaments are fundamentally different from sports contests.

Let us call this the {\em prediction tournament paradox}.
Once observed, the explanation will be quite simple.  
Possible implications for real-world prediction tournaments will be discussed in section \ref{sec;disc}.

  \section{Who wins the tournament?}
\label{sec:whowins}
Simulations in this section  use the following model for a tournament.
  \begin{quote}
  ({\bf Default tournament model}).
 There are 100 questions, with true probabilities 0.05, 0.15, 0.25, \ldots , 0.95, each appearing 10 times.
 \end{quote}
  This number of questions roughly matches the real tournaments we are aware of.  
  
  \subsection{Intrinsic variability}

\begin{figure}[h!]
  \caption{Chance variation in tournament score.}
  \label{Fig1}
  \begin{center}
  \includegraphics[width=3.0in]{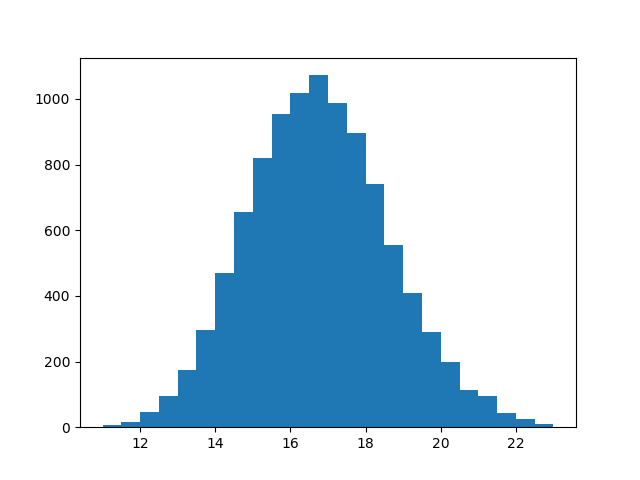}
  \end{center}
  
  \end{figure}

  In this model,  for a player who always predicted the true probabilities their mean score would be 16.75.
  But there is noticeable random variation between realizations of the tournament events, 
  illustrated in the Figure \ref{Fig1} histogram.

  The variability in  Figure \ref{Fig1} can be regarded (very roughly) as the ``luck" in the 3-part decomposition (\ref{3-part}).

 \subsection{Comparing two contestants}
We do not expect contestants to predict exactly the true probabilities, so to understand a real tournament 
we need to model inaccuracy of predictions.  This is conceptually challenging.
The basic formula (\ref{Sdef}) shows that it is the  MSE  $\sigma^2$ in forecasting  which affects score, 
so we parametrize ``inaccuracy" by the  RMS error $\sigma$.
Amongst many possible models, we take what is perhaps the simplest. 
\begin{quote}
({\bf Simple model for predictions by contestant with RMS error $\sigma$}).
When the true probability is $p$, the contestant predicts $p \pm \sigma$, each with equal probability 
(independent for different questions, and truncated to $[0,1]$). 
\end{quote}

Figure \ref{Fig2} shows the probability that, in this model tournament, a more accurate forecaster gets a better score than a  
less accurate forecaster. 
The simulation results here correspond well to intuition; indeed the probability depends, roughly, on the difference in RMS errors.

  \begin{figure}[h!]
  \caption{Chance of more accurate forecaster beating less accurate forecaster in 100-question tournament. }
  \label{Fig2}

    \begin{center}
  \begin{tabular}{|cl|llllll|}
  \hline
 &  \multicolumn{6}{ c  }{RMS error (less accurate)} & \\
 & &0.05 & 0.1 & 0.15&0.2&0.25&0.3\\
  \hline
&0&0.73&0.87&0.95&0.99&1.00&1.00\\
RMS &0.05&&0.77&0.92&0.97&0.99&1.00\\
error & 0.1&&&0.78&0.92&0.97&0.99\\
(more & 0.15&&&&0.76&0.92&0.97\\
(accurate)&0.2&&&&&0.76&0.91\\
&0.25&&&&&&0.73\\
\hline
\end{tabular}

    \end{center}
  \end{figure}

  \subsection{Rank of tournament winner}

We now consider a tournament with 300 contestants, keeping the model above for questions and forecasting accuracy.
If all contestants had equal forecasting ability  then each would be equally likely to be the winner. 
Modeling variability of accuracy amongst the field of contestants is also difficult, and again we take a simple model. 
\begin{quote}
({\bf Simple model for variability of accuracy amongst contestants}).
Abilities (measured by RMS error $\sigma$) range evenly 
across an interval, which we arbitrarily take to have length $0.3$. 
\end{quote} 
In pseudo-Python code
{\small 
\begin{quote}
{\tt

$n_c = 300$ \# number of contestants

$n_q = 100$ \# number of events

for $i = 0$ to $n_q -1$: \ $p(i) = 0.05 + 0.1 \lfloor i/10 \rfloor$ \# True probability of event $i$

\hspace*{0.3in} $B(i) = $  Bernoulli($p(i)$): \ \# outcome of event $i$

\# $[\sigma_0, \sigma_0 + 0.3]$  range of RMS errors for different contestants

for $j = 0$ to $n_c - 1: \  \sigma(j) = \sigma_0 + 0.3 j/n_c$: \ \# RMS error of contestant $j$

\hspace*{0.3in} for $i = 0$ to $n_q -1$: \  $q(i,j) = p(i) \pm \sigma(j)$  \#prediction of contestant $j$ \\     \hspace*{3.3in}    for event $i$,   with random  $\pm$ 

\hspace*{0.3in}  \hspace*{0.3in} score$(i,j) =  (B(i) - q(i,j))^2$ \#squared error, contestant $j$ event $i$

\hspace*{0.3in} score$(j) = \sum_i$ score$(i,j)$ \# total score of contestant $j$.
}
\end{quote}
}

\medskip
\noindent
With 300 contestants, the top-ranked ability is little different from that of the second- or third-ranked, so the chance of  the top-ranked contestant winning will not be large in absolute terms.
But common sense and the Figure \ref{Fig2} results suggest the winner will be one of the relatively top-ranked contestants; as in any sport, the probability of being the tournament winner should decrease with rank of ability.

\begin{figure}[h!]
  \caption{Rank of tournament winner, 300 contestants, error parameters $0 < \sigma < 0.3$ }
  \label{Fig3}
  
  \begin{center}
  \includegraphics[width=3.0in]{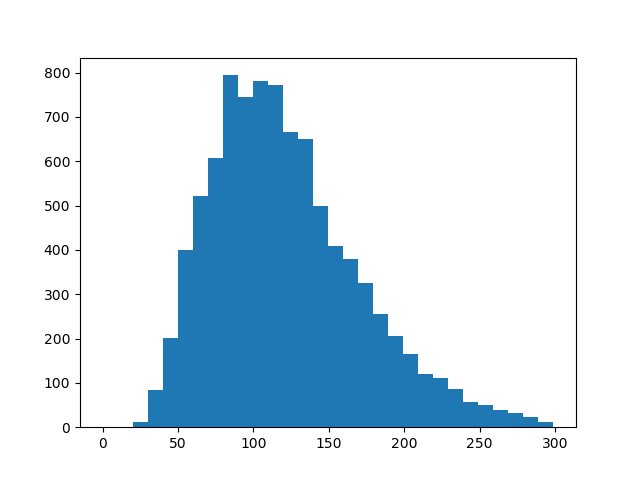}
  \end{center}
  
  \end{figure}

Figure \ref{Fig3} shows the  results of the first simulation we did, taking the interval of RMS error parameters to be $[0,0.3]$.

\begin{figure}[h!]
  \caption{Rank of tournament winner, 300 contestants, error parameters $0.05 < \sigma < 0.35$ (left) and  $0.1 < \sigma < 0.4$ (center)  and $0.15 < \sigma < 0.45$ (right). }
  \label{Fig4}
  
  \begin{center}
  \includegraphics[width=1.9in]{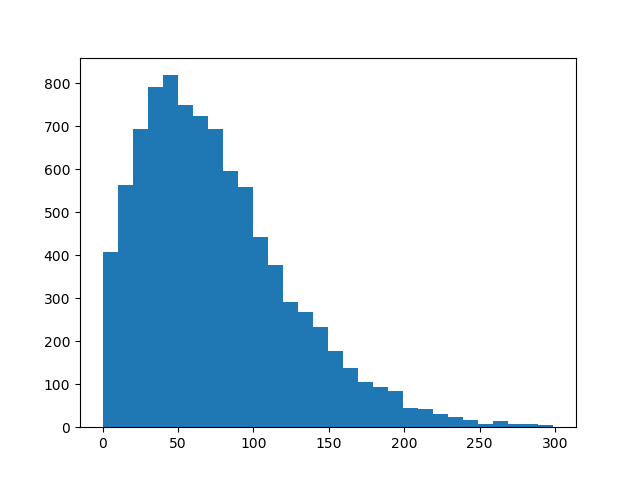}
  \includegraphics[width=1.9in]{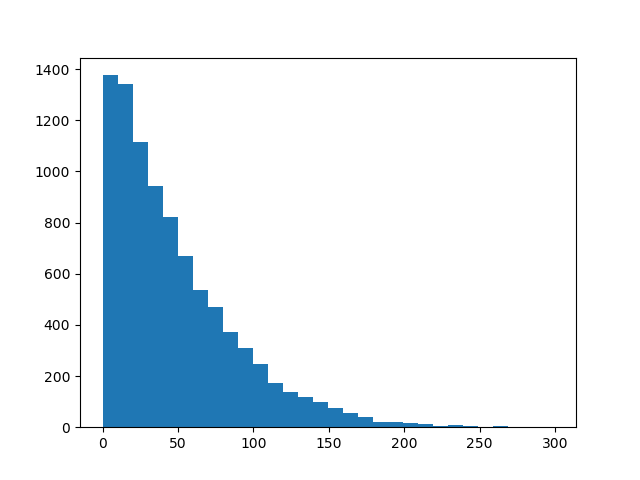}
   \includegraphics[width=1.9in]{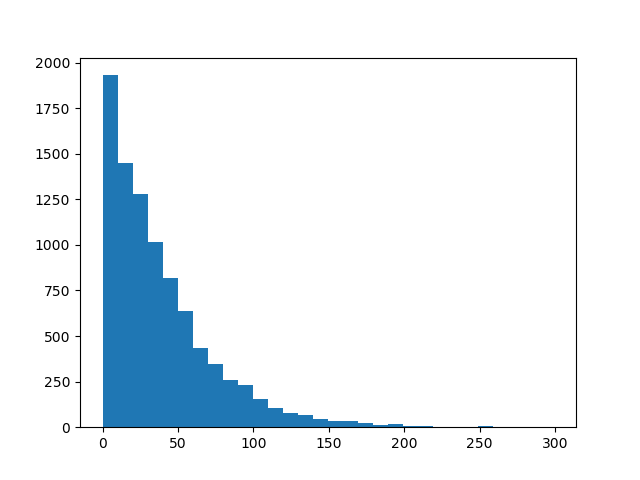}
  \end{center}
  
  \end{figure}

So here the winner is relatively most likely to be around the 100th most accurate of the 300 contestants, and the top-ranked contestants never win.  
 This is in striking contrast to intuition -- a paradox, in that sense.
 Indeed one might well suspect an error in coding the simulations.   
 However if we shift the assumed interval of $\sigma$  successively to 
 $[0.05, 0.35]$ and $[0.1,0.4]$ and $[0.15,0.45]$ then 
 (see Figure \ref{Fig4}) we do soon see the intuitive ``winning probability decreases with rank"  property, but still the winners are not as strongly concentrated amongst the very best forecasters  as one might have guessed.
  
  \paragraph{First explanation of the paradox.}
  Once observed, the original paradox is easy to explain in words.  
  In the specific setting of  Figure \ref{Fig3} the handful of top-rated contestants are making almost exactly the same predictions and therefore getting almost exactly the same score -- as if there were just one such contestant.
But looking at contestants with $\sigma$ around 0.1 they are making slightly different predictions, on average scoring less well;
 but by chance, 
for some contestants, most of the predictions will vary in the direction of the outcome that actually occurred, and so these contestants will get a better score by pure luck.
As a physical analogy, imagine contestants who each shoot successively at 100  different red targets. 
But there is an invisible-to-contestants blue target randomly displaced from each red target, and they are scored by the average distance between the shot and the blue target.
The skillful contestants whose shots land close to the red targets will all get roughly the same score.  
Less skillful contestants will typically get lower scores, but some will by chance have more errors in the directions toward the blue targets and therefore, by luck, get a better score.
This is an instance of a  mean-variance trade-off.  
See section \ref{sec;disc} for further discussion.
  
\begin{figure}[h!]
  \caption{Alternate model: rank of tournament winner, 300 contestants, error parameters $0.05 < \sigma < 0.35$ (left) and  $0.1 < \sigma < 0.4$ (center)  and $0.15 < \sigma < 0.45$ (right). }
  \label{Fig4.5}
  
  \begin{center}
  \includegraphics[width=1.9in]{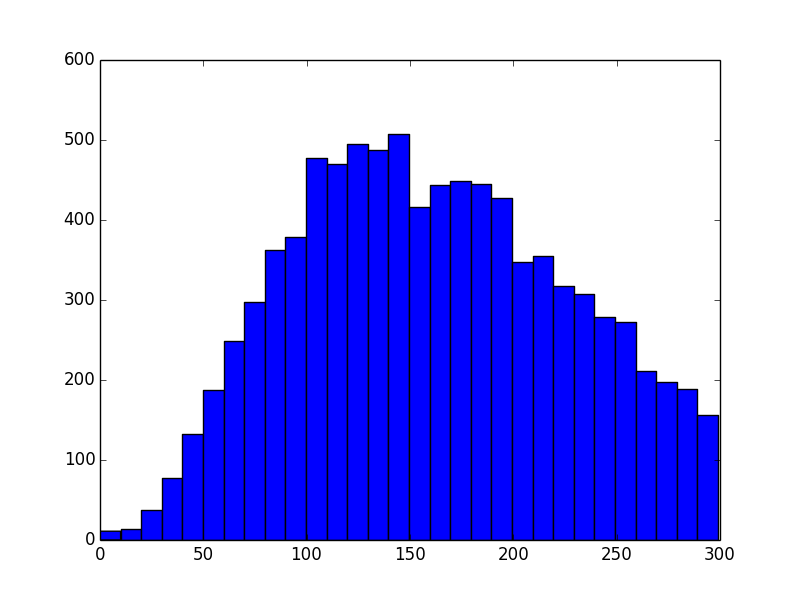}
  \includegraphics[width=1.9in]{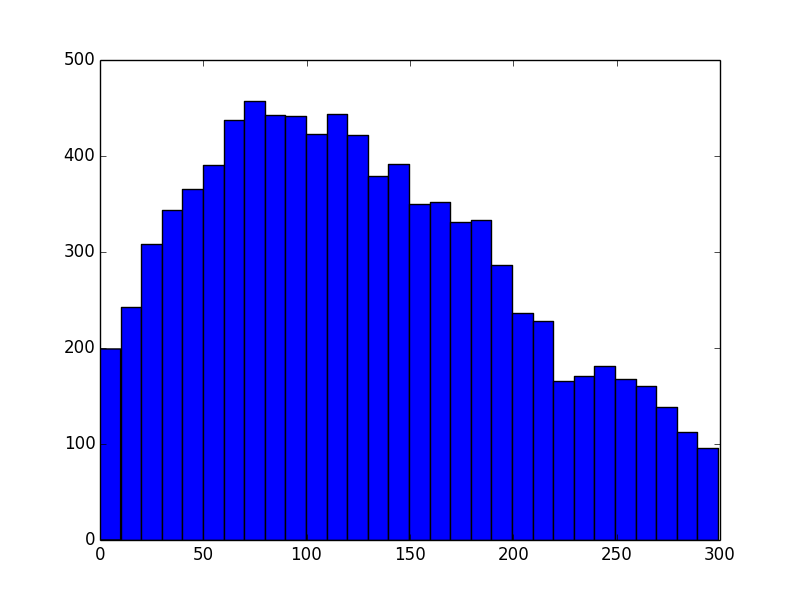}
   \includegraphics[width=1.9in]{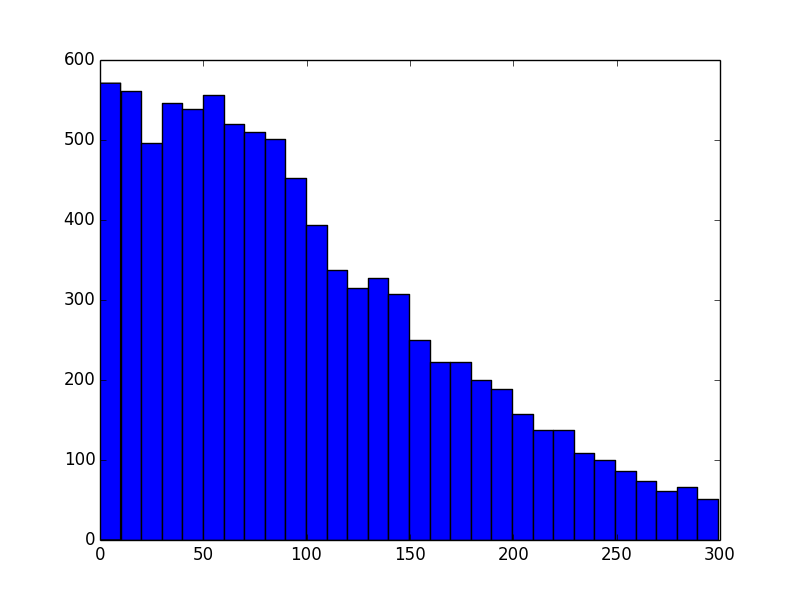}
  \end{center}
  
  \end{figure}

  One might wonder whether the specific implementation of RMS error $\sigma$ as `` predict $p \pm \sigma$ with equal probability" 
  had any effect.  Simulations with  the alternative ``predict a random variable uniform in $[p- \sigma \sqrt{3}, p+ \sigma \sqrt{3}]$" implementation
  are shown in Figure \ref{Fig4.5}, and the effect is even stronger in that model.

\begin{figure}[h!]
  \caption{Ranks of tournament winner (left) and  
  of top ten finishers (right), $0.15 < \sigma < 0.45$. }
  \label{Fig5}
  
  \begin{center}
   \includegraphics[width=2.2in]{figure_spread_0,15-0,45.png}
      \includegraphics[width=2.2in]{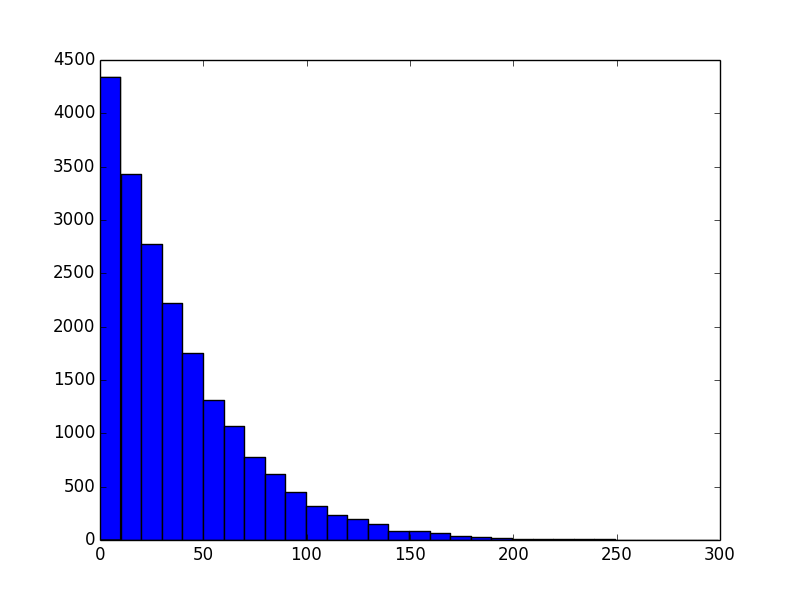}
  \end{center}
  
  \end{figure}

  One might also wonder if this behavior is special to the winner, but we see a similar effect if we look at the ranks of contestants whose scores are in the top 10 -- see Figure \ref{Fig5}.

  \newpage

\begin{figure}[h!]
  \caption{Rank of tournament winner, $0 < \sigma < 0.3$ (left) and  $0.15 < \sigma < 0.45$  (right) for systematic over- or under-estimation. }
  \label{Fig6}
  
  \begin{center}
  \includegraphics[width=2.9in]{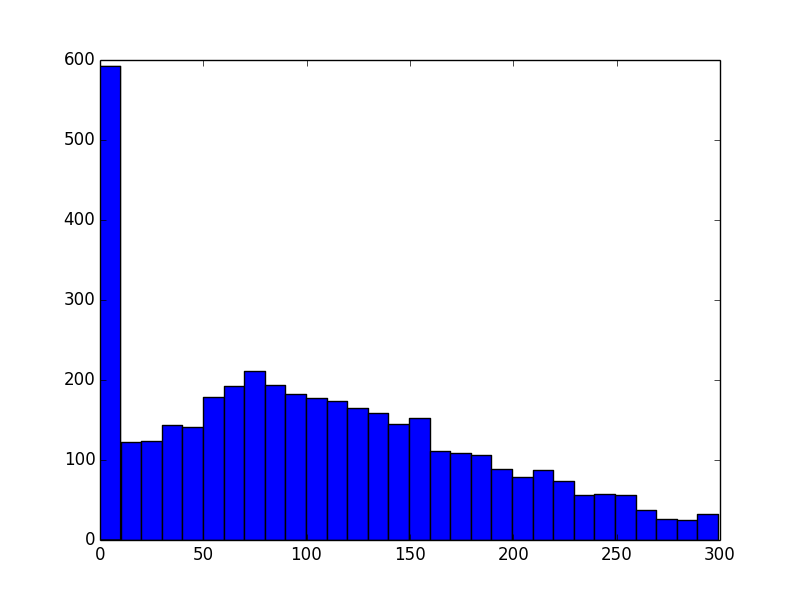}
    \includegraphics[width=2.9in]{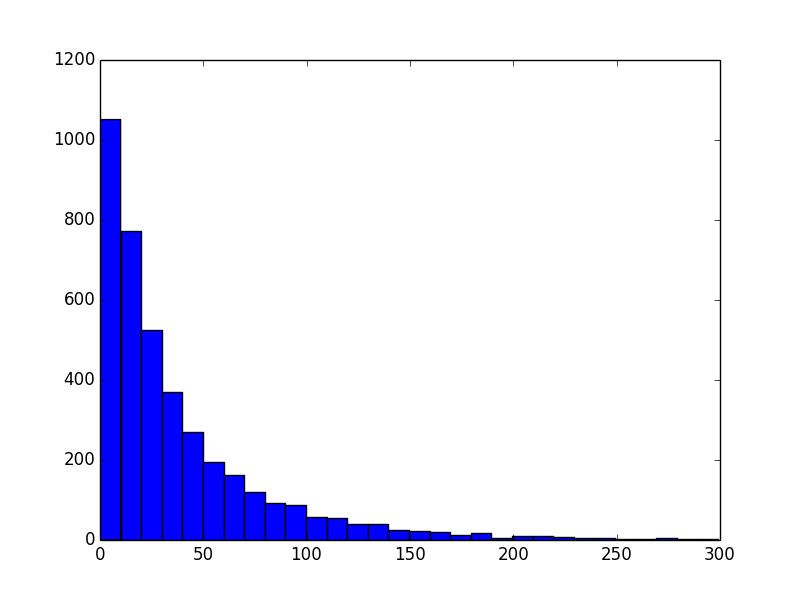}
  \end{center}
  
  \end{figure}

Are these results merely artifacts of the specific model?
Suppose instead we took a ``one-sided" model of inaccuracy,  say in the sense that contestants systematically over-estimate probabilities.
Then there would be a comparatively large chance that 
the top-ranked contestant is the winner, from the case that more event outcomes than expected turned out to be ``no". 
What about a model in which  half the contestants systematically over-estimate probabilities and the other half
systematically under-estimate? 
The results are shown in Figure  \ref{Fig6}.  
As before the inaccuracy parameter $\sigma$ of different contestants varies evenly over $[0,0.3]$, 
but now the prediction model is

\begin{quote}
for half the contestants,
when the true probability is $p$, the contestant predicts a random value uniform in $[p, p+ \sigma \sqrt{3}]$ 
(independent for different questions, and truncated to $[0,1]$);
for the other half of the contestants,
when the true probability is $p$, the contestant predicts a ran.om value uniform in $[p- \sigma \sqrt{3},p]$ 
(independent for different questions, and truncated to $[0,1]$);
\end{quote}

Here, in the case $0 < \sigma < 0.3$ with accurate forecasters, we see a combination of the effects noted above. A near-top-rated contestant will likely win when 
the pattern of event outcomes is relatively close to balanced 
(events of probability $p$ happen a proportion $p$ of the time), but as in the previous one-sided case 
some of the biased contestants will, by luck, do better when outcomes are unbalanced. 
In the case $0.15 < \sigma < 0.45$ of all inaccurate forecasters the pattern of winner ranks is similar to the original model (Figure \ref{Fig4}).

  \section{Discussion}
  \label{sec;disc}
   One subtle point not discussed earlier is that our measure $\sigma$ of accuracy is implicitly a long-term average.
 By modeling predictions as random, each  contestant in our model has an empirical RMS prediction error in a finite tournament.
 In the size of tournaments simulated here, for contestants with equal $\sigma$ 
the correlation between score and RMS empirical prediction error is small: the ``noise" of event outcomes is overwhelming.
One might hope that differences in long-term accuracy would show up in reasonable size tournaments, 
but part of our ``paradox" is the observation  that 100-question tournaments are not sufficient.

\paragraph{Explanations of the paradox.}
To elaborate our earlier mean-variance-tradeoff explanation of the paradox 
with some numbers, consider a 100-question tournament in which the true probabilities
are all 0.5.  So a perfectly accurate forecaster will score exactly 25.0.
Now consider a contestant who predicts 0.4 or 0.6 randomly on each question.
Their score is random with expectation 26 and standard deviation 0.98, so have around 15\% chance
to beat the perfect forecaster.  
If instead predictions were 0.3 or 0.7, the expectation and s.d. become  29 and 1.83.
Moreover, as a special feature of the ``all true probabilities are 0.5" setting, 
different contestants' scores are independent.  
In our simulated  setting of 300 contestants  with RMS prediction errors ranging from 0 to 0.3.
some scores will by chance be around 3 s.d.'s below expectation, and by this  
back-of-an-envelope argument we expect a winning score around 23
and we will not be surprised if this comes from  the 100th or 200th best forecaster.

Our simulations used the more plausible {\em default probability model} with varying true probabilities, 
and here there is a specific complicated dependence structure for the scores of different contestants, not amenable to 
convincing back-of-an-envelope calculations or convincing asymptotic approximations or human-interpretable algebra,
which is why we have relied on simulations. 

  Experts in statistical methodology might readily think of their own explanations of
  or analogies for the paradox.  
For instance it is partly analogous to standard 
  multiple comparison settings, though the specific dependence structure arising in this  
``estimating probabilities" context
is quite different from the usual contexts 
of multiple comparisons for experimental or observational data.

In a typical sports setting, the winner of a tournament is indeed relatively more likely to be one of the best teams.
So it is important to realize how our prediction tournament setting 
is conceptually very different 
from the more familiar setting of a contest in which each contestant earns points (directly reflecting skill, 
as in a basketball shot) in each of 100 rounds and the winner is the contestant with the most total points.
In sports an ``error" -- that is, not making the percentage play -- is usually costly and only rarely is it luckily beneficial 
(a soccer shot that would miss the goal might luckily be deflected by a defender into the goal, for instance).
But in the probability prediction context, predicting a 60\% or 40\% probability when the true probability is 50\% is almost 
equally beneficial or costly.   
Loosely speaking, errors in predicting probabilities have only a second-order effect:  
100 errors in a sequence of sports matches are more costly than 100 errors in predicting probabilities, 
and the latter might indeed by pure luck be overall beneficial.

   \paragraph{Practical relevance?} 
  Our model is over-simplified in many ways; does it have implications for real-world prediction tournaments?
  Currently (announced February 2018) IARPA is offering \$200,000 in prizes for top performers in its upcoming 
  Geopolitical Forecasting Challenge (\cite{IARPA-announce}); no doubt this will encourage volunteers to participate, but is it effective in identifying the best forecasters?
  
 The authors of \cite{essay}  write ``some forecasters are, surprisingly consistently, better than others", 
 and background to this assertion can be found in  \cite{identifying}:
  \begin{quote}
 [the winning strategy for teams over several successive tournaments was] culling off top performers each year and assigning them into elite teams of superforecasters. Defying expectations of regression toward the mean 2 years in a row, superforecasters maintained high accuracy across hundreds of questions and a wide array of topics.
 \end{quote}

Designers of that strategy were implicitly assuming that doing well in a tournament is strong evidence for ability (rather than luck), though our model results  suggest that this assumption 
  deserves some scrutiny. 
  However a main focus of the recent literature  is arguing for the effectiveness of training methods, so that (if it were correct to downplay the effectiveness of ``culling off top performers each year" in selecting for prior ability) 
  our results actually reinforce that argument.

 A superficial conclusion of our results is that winning a prediction tournament is strong evidence of superior ability 
  {\em only} when the better forecasters' predictions are {\em not} reliably close to the true probabilities.\footnote{Asking whether ``close to the true probabilities"  is true in practice leads to basic issues in the philosophy of the meaning of {\em true probabilities}, not addressed here.}
  But are our models realistic enough to be meaningful?
  Two features of our ``simple model for predictions by contestant with RMS error $\sigma$" are unrealistic.
  One is that contestants have no systematic bias towards too-high or too-low forecasts.
  A more serious issue is that the errors are assumed independent over both questions and contestants.  
  In reality, if all contestants are making judgments on the same evidence, then (to the extent that relevant evidence is incompletely known) there is surely a tendency for most contestants to be biased in the same direction on any given question. 
  Implicit in 
  our model (and in our ``first explanation" previously) is that, in a large tournament, this ``independence of errors" assumption means that different contestants will explore somewhat uniformly over the space of possible prediction sequences close to the true probabilities, 
  whereas in reality one imagines the deviations would be highly non-uniform.

  For recent relevant technical literature see \cite{wilkowski17,wilkowski} and citations therein.  
 In particular, when viewed as game theory with each player's only objective being to win the tournament, under the usual scoring scheme the optimal strategy involves {\em not} making truthful predictions, so one can study alternative scoring schemes that incentivize truthful reporting and are more likely to identify the best forecasters.



\end{document}